\documentclass[11pt]{article}
\usepackage[margin=1in]{geometry}
\usepackage{amsmath,amssymb,amsthm,booktabs}
\usepackage[hidelinks]{hyperref}
\hypersetup{pdftitle={A bordism theory related to matrix Grassmannians},pdfauthor={A. V. Ershov}}
\newtheorem{theorem}{Theorem}[section]
\newtheorem{proposition}[theorem]{Proposition}
\newtheorem{lemma}[theorem]{Lemma}

\theoremstyle{definition}
\newtheorem{example}[theorem]{Example}
\newtheorem{remark}[theorem]{Remark}
\newcommand{\Gr}{\mathrm{Gr}}
\newcommand{\SU}{\mathrm{SU}}
\newcommand{\BSU}{B\SU}
\newcommand{\fr}{\mathrm{fr}}
\newcommand{\kk}{K^0}
\newcommand{\CC}{\mathbb C}
\newcommand{\ZZ}{\mathbb Z}
\newcommand{\QQ}{\mathbb Q}
\newcommand{\RR}{\mathbb R}
\newcommand{\triv}[1]{\underline{#1}}
\newcommand{\units}{\mathcal U}
\DeclareMathOperator{\rank}{rank}
\DeclareMathOperator{\chernchar}{ch}
\DeclareMathOperator{\Th}{Th}
\DeclareMathOperator{\Tel}{Tel}
\title{A bordism theory related to matrix Grassmannians}
\author{A.\,V. Ershov}
\date{Revised manuscript}
\begin{document}
\maketitle
\begin{abstract}
Matrix Grassmannians give finite-dimensional representatives for virtual bundles of rank one and vanishing first Chern class. We use them to study two groups of labelled manifolds. A choice of stable normal framing gives framed bordism of $\BSU_{\otimes}$ and an integral Pontryagin--Thom description. Forgetting the framing gives a quotient defined by the boundary relations for labelled, stably parallelizable manifolds. The forgetful map is surjective with torsion kernel. Both rational bordism rings are polynomial on explicitly labelled spheres, one generator in each even dimension at least four. We compute the characteristic numbers of their products and show that each sphere generator spans an integral direct summand in its own degree.
\end{abstract}

\section{Introduction}
A unital copy of $M_k(\CC)$ in $M_{kl}(\CC)$ determines a second matrix algebra, its centralizer. When $k$ and $l$ are relatively prime, a family of such subalgebras carries two special unitary vector bundles $V$ and $W$ whose tensor product is canonically trivial. In a range controlled by $\min(k,l)$, the family represents a single virtual bundle $\xi$ of rank one. The two factors encode $\xi$ and its multiplicative inverse by the identities
\[
[V]=k\xi,\qquad [W]=l\xi^{-1}.
\]
This interpretation of the matrix Grassmannian was developed in \cite{Er2000,Er2003}. It allows the tensor product of labels to be seen on finite-dimensional algebras. The bordism question asks which labelled manifolds can occur as boundaries and what information is detected by the label's characteristic numbers. It belongs to the study of bordism of bundles exemplified by Floyd \cite{Floyd}.

There are two ways a stable trivialization enters this question. A choice of normal framing produces a Pontryagin--Thom collapse. The existence of a normal framing is also a property one can require of the underlying manifold without recording a choice. In the latter case we impose the boundary relations as relations in an abelian group. The difference between the two groups is already visible on an oriented circle with the constant label. It disappears after passage to rational coefficients.

Write $B=\BSU_{\otimes}$ for the space $\BSU$ with multiplication induced by tensor product of rank-one virtual $\SU$-bundles. Write $\Omega^{\fr}_*(B)$ for framed bordism with maps to $B$ and $A_*(B)$ for the group with forgotten framings defined in Section~\ref{sec:bord}. Our main conclusions are as follows.
\begin{enumerate}
\item If $k,l\ge2$ are coprime, $d<2\min(k,l)$ and $N=2kl$, then there are integral isomorphisms
\[
  \pi_{d+N}\!\left(\Th(\vartheta_{k,l})\right)
  \cong\Omega_d^{\fr}(B)\cong\pi_d(\Sigma^\infty B_+),
\]
where $\vartheta_{k,l}$ is the tautological tensor product on $\Gr_{k,l}$ and $B_+$ denotes the space with a disjoint basepoint.
\item Forgetting a framing induces a surjective ring homomorphism $\Omega^{\fr}_*(B)\to A_*(B)$ with torsion kernel. In particular the two integral rings have isomorphic quotients by their torsion subgroups.
\item After tensoring with $\QQ$, either ring is $\QQ[t_2,t_3,\ldots]$, where $t_n$ has degree $2n$ and is represented by $S^{2n}$ with label $1+b_n$ for a normalized Bott generator $b_n$. In both integral groups $t_n$ spans a direct summand isomorphic to $\ZZ$.
\end{enumerate}

The geometric cycles and the rational generators originate in the 2010 preprint \cite{Er2010}. The framing in the first assertion and the definition by generated boundary relations in the second give precise versions of its two geometric constructions. The algebraic normalization in Section~\ref{sec:matrix} also specifies the effect of a homotopy of trivializations on a matrix Grassmannian representative.

\section{Rank-one labels and finite matrix Grassmannians}\label{sec:matrix}
For a finite CW complex $X$ put
\[
 \units(X)=\{\xi\in\kk(X):\rank(\xi)=1,\ c_1(\xi)=0\},
\]
where ranks are taken on each component. The rank-zero ideal is nilpotent by the skeletal filtration. Thus every rank-one class is a unit, and the first Chern class formula for tensor products makes $\units(X)$ an abelian group. Adding a trivial line to a reduced stable $\SU$-class identifies $[X,B]$ with $\units(X)$. Here one may use the determinant map $BU\to BU(1)$, its line-bundle section and stable Whitney subtraction of the determinant line to identify $[X,\BSU]$ with the reduced classes of first Chern class zero. The product on $B$ corresponds to multiplication in $\kk(X)$. For labels on distinct bases we use the external tensor product $\boxtimes$.

Fix relatively prime integers $k,l\ge2$. Let $G_{k,l}=\Gr_{k,l}$ denote the space of unital $*$-subalgebras of $M_{kl}(\CC)$ isomorphic to $M_k(\CC)$. Unitary conjugation is transitive. The stabilizer in the projective unitary group is the product of the two projective unitary groups acting on the two factors. In $\SU(kl)$ the tensor homomorphism $\SU(k)\times\SU(l)\to\SU(kl)$ has trivial kernel, since $k$ and $l$ are coprime. Its image contains the full center of $\SU(kl)$, since the scalar subgroups of orders $k$ and $l$ generate that center. Consequently
\begin{equation}\label{eq:homogeneous}
 G_{k,l}=\SU(kl)/(\SU(k)\otimes\SU(l)).
\end{equation}
The two standard representations of the stabilizer give bundles $V_{k,l}$ and $W_{k,l}$. The tensor representation extends to $\SU(kl)$, giving a specified special unitary trivialization
\begin{equation}\label{eq:univtriv}
 \vartheta_{k,l}=V_{k,l}\otimes W_{k,l}
 \cong G_{k,l}\times\CC^{kl}.
\end{equation}
For example, in the associated bundle description it is $[g,z]\mapsto (gH,gz)$, where $H=\SU(k)\otimes\SU(l)$. The endomorphism algebras of $V_{k,l}$ and $W_{k,l}$ are the tautological subalgebra and its centralizer.

A map from $X$ to $G_{k,l}$ amounts to a reduction of the trivial principal $\SU(kl)$-bundle to $H$. Equivalently it specifies a pair of $\SU$-bundles $(V,W)$ and a determinant-preserving unitary trivialization of $V\otimes W$, up to factor isomorphisms carrying that trivialization to the given one. A homotopy of maps amounts to such data on $X\times[0,1]$ with the prescribed endpoint identifications. This description also follows from the homotopy fiber sequence
\begin{equation}\label{eq:fibration}
 G_{k,l}\longrightarrow B\SU(k)\times B\SU(l)
                \xrightarrow{\ \otimes\ }B\SU(kl).
\end{equation}

\begin{lemma}[Recovering the label]\label{lem:normalize}
Let $V,W$ be $\SU$-bundles of coprime ranks $k,l$ over a finite CW complex $X$, with $[V][W]=kl$ in $\kk(X)$. There is a unique $\xi\in\units(X)$ satisfying $k\xi=[V]$ and $l\xi^{-1}=[W]$. It is natural under pullback.
\end{lemma}
\begin{proof}
Put $R=\kk(X)$. In $R[1/k]$ take $x=[V]/k$, and in $R[1/l]$ take $y=([W]/l)^{-1}$. The inverse exists because the rank-zero ideal is nilpotent. The two classes agree in $R[1/(kl)]$ by the product identity. For coprime $k,l$, any abelian group $L$ is the equalizer of its two localizations $L[1/k],L[1/l]$ over $L[1/(kl)]$. Indeed, an element mapping to zero in both is killed by coprime powers of $k$ and $l$. For existence write matching elements as $a/k^r,b/l^s$. Increasing both exponents and multiplying the numerators accordingly gives $l^s a=k^r b$ in $L$. If $u k^r+v l^s=1$, then $u a+v b$ maps to both elements.

Gluing $x$ and $y$ gives a unique rank-one class $\xi$. Both required identities hold in each localization, hence in $R$. Since $V,W$ have vanishing first Chern class, $k c_1(\xi)=0$ and $l c_1(\xi^{-1})=-l c_1(\xi)=0$. Coprimality gives $c_1(\xi)=0$. The two identities prescribe $x$ and $y$, which also proves uniqueness and naturality.
\end{proof}

Apply Lemma~\ref{lem:normalize} on the compact homogeneous space $G_{k,l}$ to obtain a universal class $u_{k,l}$. Choose a representing map $j_{k,l}:G_{k,l}\to B$, unique up to homotopy. Applied to the standard subalgebra at a point, this class is the unit label.

\begin{theorem}[Finite-dimensional representatives]\label{thm:range}
The comparison $j_{k,l}$ induces isomorphisms on $\pi_i$ for $i\le 2\min(k,l)$. In particular, for a finite CW complex $X$ with $\dim X\le 2\min(k,l)$, it gives a natural bijection
\[
 [X,G_{k,l}]\ \xrightarrow{\ \cong\ }\ \units(X),\qquad
 [f]\longmapsto f^*u_{k,l}.
\]
The factors associated to $f$ represent $k f^*u_{k,l}$ and $l(f^*u_{k,l})^{-1}$.
\end{theorem}
\begin{proof}
Both spaces are simply connected. The quotient in \eqref{eq:homogeneous} also has $\pi_2=0$. In positive odd degrees $s\ge3$ with $s<2\min(k,l)$, the map of stabilizer to ambient group on $\pi_s$ is
\[
 \ZZ^2\to\ZZ,\qquad (a,b)\mapsto la+kb.
\]
Indeed, on each factor tensoring by the identity is a block sum of $l$ or $k$ copies. The indicated groups and the vanishing of the even groups in the same range follow from Bott periodicity and $\SU(r+1)/\SU(r)=S^{2r+1}$ \cite{Atiyah}. The displayed map is surjective. The exact sequence of \eqref{eq:homogeneous} now gives vanishing in odd degrees in the asserted range.

In degree $i=2m\ge4$ with $i\le2\min(k,l)$, it identifies $\pi_i(G_{k,l})$ with the subgroup $\ZZ(k,-l)$ in the two groups classifying the factors on $S^i$. Here only the factor groups in degree $i-1$ have to be stable. The ambient group satisfies $\pi_i\SU(kl)=0$ because $i<2kl$. On a generator the first factor has reduced $K$-class $k$ times a Bott generator. As $k u_{k,l}=[V_{k,l}]$, the comparison takes the chosen generator to one Bott generator of $\pi_i B$. This proves the assertion including the even endpoint. For $e=2\min(k,l)$ the induced map of $e$th Postnikov sections is a homotopy equivalence. A complex of dimension at most $e$ has the same homotopy classes of maps into a space and into its $e$th Postnikov section, by obstruction theory. This proves the bijection. The last assertion is Lemma~\ref{lem:normalize}.
\end{proof}

For a geometric interpretation, recall that $B\SU(r)\to\BSU$ is $(2r+1)$-connected for $r\ge2$, by the same unitary stability calculation. Thus, in the range of Theorem~\ref{thm:range}, the rank-$k$ and rank-$l$ factors are individually determined up to special unitary isomorphism by their stable classes. Conversely, given two such factors whose tensor product is trivial as an $\SU(kl)$-bundle, each determinant-preserving trivialization defines a map to $G_{k,l}$. All these maps have the same label by Lemma~\ref{lem:normalize}, and hence are homotopic by Theorem~\ref{thm:range}. This is a statement about the homotopy classes of the resulting maps.

\begin{remark}[Trivializations of fixed factors]
Two trivializations can determine the same matrix Grassmannian class. For the trivial factors of ranks $2$ and $3$ on $S^3$, fixed-factor special unitary trivializations have components $\pi_3\SU(6)=\ZZ$. Changes of coordinates of the factors act by the image of $(a,b)\mapsto3a+2b$, which is all of $\ZZ$. Accordingly $[S^3,G_{2,3}]$ is a singleton. This accounts for the uniqueness asserted in Theorem~\ref{thm:range}.
\end{remark}

\paragraph{Products and stabilization.}
For $\gcd(km,ln)=1$ the tensor product of subalgebras gives
\[
 G_{k,l}\times G_{m,n}\longrightarrow G_{km,ln},\qquad (P,Q)\longmapsto P\otimes Q.
\]
We identify the ambient matrix algebras by the usual permutation of tensor factors. The two bundles pull back to the external tensor products of the corresponding factors, as follows either from the stabilizer representations or directly from the reductions. Lemma~\ref{lem:normalize} then identifies the pullback of $u_{km,ln}$ with $u_{k,l}\boxtimes u_{m,n}$. Fixing $Q=M_m(\CC)\otimes1_n$ gives the stabilization. A convenient tower is $G_r=G_{2^r,3^r}$ for $r\ge1$, with stabilization by the standard point of $G_1$. Products take values in $G_{r+s}$. The compatibility of the universal labels supplies homotopies between successive comparison maps. They give a map $\Tel_r G_r\to B$ which is a weak equivalence by Theorem~\ref{thm:range}, since homotopy groups of the telescope are the direct limits of the homotopy groups of the stages. This tower supplies all finite-dimensional representatives needed below.

\section{Stable framings and the boundary relations}\label{sec:bord}
Let $Y$ be a CW space. A cycle for $\Omega_d^{\fr}(Y)$ is a closed smooth $d$-manifold $M$, a stable normal framing $\mathfrak f$, and a map $f:M\to Y$. A framing is a homotopy class of trivializations after adjoining trivial real summands to a normal bundle for an embedding in Euclidean space. It determines an orientation of $M$. A bordism carries a stable normal framing inducing the assigned end framings with the usual boundary convention, and a map to $Y$ extending the end maps. Equivalently one can use neat embeddings and product collars at the boundary. For $Y=B$ the product takes the Cartesian product with the product framing and the external tensor product of the labels.

Define $A_d(Y)$ to be the abelian group on orientation-preserving diffeomorphism classes of pairs $(M,f)$ where $M$ is closed, smooth, oriented and stably parallelizable, and $f:M\to Y$ is taken up to homotopy. Impose the relations
\begin{align*}
 [M\sqcup M',f\sqcup f']&=[M,f]+[M',f'],\\
 [\partial W,F|_{\partial W}]&=0
\end{align*}
for every compact smooth oriented stably parallelizable $W^{d+1}$ and map $F:W\to Y$. Stable parallelizability is equivalent to stable normal triviality. The cylinder relation identifies $[-M,f]$ with $-[M,f]$. For $Y=B$ products descend because $W\times M'$ is stably parallelizable with boundary $\partial W\times M'$ and the external product of the labels extends over it. This gives a graded-commutative ring with unit the positively oriented point with unit label.

The use of generated boundary relations makes the definition independent of a gluing assertion for arbitrary unframed bordisms. For example, the nontrivial oriented $S^2$-bundle over $S^2$ is glued from two stably parallelizable pieces $D^2\times S^2$. To see the obstruction, realize the bundle as $p:S(E)\to S^2$ for an oriented real $3$-plane bundle $E$ with nonzero $w_2$. Such a bundle is obtained as a complex line of Chern number one, regarded as a real plane, plus a trivial real line. It has a section of its sphere bundle, so $p^*$ is injective on cohomology. The isomorphism $TS(E)\oplus\triv{\RR}\cong p^*(TS^2\oplus E)$ gives $w_2(TS(E))=p^*w_2(E)\ne0$. Compatible framings on collars ensure the usual framed gluing theorem. The example concerns a specified gluing and leaves open the question of replacing a chain of unframed bordisms by a single one with the same ends.

\begin{theorem}[Integral Thom model]\label{thm:thom}
For a CW space $Y$ there is a natural isomorphism
\[
 \Omega_d^{\fr}(Y)\cong\pi_d(\Sigma^\infty Y_+).
\]
For $d\ge0$, $d<2\min(k,l)$ and $N=2kl$, this gives isomorphisms
\begin{equation}\label{eq:finitePT}
 \pi_{d+N}(\Th(\vartheta_{k,l}))\ \cong\
 \Omega_d^{\fr}(G_{k,l})\ \xrightarrow[\cong]{(j_{k,l})_*}\ \Omega_d^{\fr}(B).
\end{equation}
They are compatible with stabilization of the matrix factors. The stable product is induced by the multiplication of $B$.
\end{theorem}
\begin{proof}
The first assertion is the labelled Pontryagin--Thom theorem \cite{Cohen}. Briefly, a framed embedding in $\RR^{d+q}$ and its label give the collapse
\[
 S^{d+q}\longrightarrow\Th(\nu_M)\xrightarrow{\mathfrak f}
 S^q\wedge M_+\xrightarrow{1\wedge f_+} S^q\wedge Y_+
\]
after stabilization. Transversality in the normal coordinates recovers the normally framed submanifold with its map to $Y$. The relative construction gives the bordism relation. Exterior products correspond to smash products of collapse maps.

By \eqref{eq:univtriv}, $\Th(\vartheta_{k,l})=S^N\wedge(G_{k,l})_+$, up to the indicated choice of standard coordinates. This space is $(N-1)$-connected. As $d+N\le2N-2$, Freudenthal identifies its homotopy group in degree $d+N$ with the stable group. Desuspending $N$ times gives the first isomorphism in \eqref{eq:finitePT}. To obtain the second, take a fibration replacement of $j_{k,l}$. Its fiber is $d$-connected by Theorem~\ref{thm:range}, since $d+1\le2\min(k,l)$. Thus maps of $d$-dimensional cycles lift up to homotopy, and maps of $(d+1)$-dimensional bordisms lift relative to specified boundary lifts. The underlying framings are carried along in both constructions.

If $N'=2klmn$, the map between the Thom groups is suspension by $N'-N$, followed by the map induced on the bases by $G_{k,l}\to G_{km,ln}$. In the canonical coordinates \eqref{eq:univtriv} this gives the map of trivial Thom spaces $S^{N'}\wedge(G_{k,l})_+\to S^{N'}\wedge(G_{km,ln})_+$. Its compatibility with \eqref{eq:finitePT} follows from the pullback identity for the labels. For two cycles the normal framings combine by direct sum, while the two labels combine by external tensor product. This proves the product assertion.
\end{proof}

\begin{example}[The circle]\label{ex:circle}
Every label on an oriented circle is homotopic to the constant label because $B$ is simply connected. Each such circle bounds a labelled disk. Consequently $A_1(B)=0$. In the framed theory the standard point of $G_{2,3}$ gives a retraction onto the Thom sphere $S^{12}$. Hence $\pi_{13}(S^{12})=\ZZ/2$ injects into $\pi_{13}(\Th(\vartheta_{2,3}))$. In fact $\Omega_1^{\fr}(B)=\ZZ/2$. The basepoint splitting of $\Sigma^\infty B_+$ has sphere summand and the reduced summand has zero first homotopy group since $B$ is simply connected. Thus forgetting the framing changes the integral group even for a constant label.
\end{example}

\begin{example}[Real and complex triviality]\label{ex:six}
A framed cycle may be represented by its framing and a map to $G_{k,l}$. It is essential that \eqref{eq:univtriv} is a \emph{complex} trivialization. To see the distinction, equip $S^6$ with its standard almost complex structure $J$ \cite{KP}. Take $V=(TS^6,J)\oplus\triv{\CC}$ and $W=\triv{\CC}^{\,5}$. These bundles admit special unitary structures and have coprime ranks $4$ and $5$. The real bundle $V_{\RR}$ is trivial since $TS^6\oplus\triv{\RR}$ is trivial. Hence the underlying real bundle of $V\otimes W$ is trivial of rank $40$, and can be identified with the normal bundle of the standard $S^6\subset\RR^{46}$. For the cohomology generator $v$ evaluating to one on the complex orientation, however,
\[
 c_3(V\otimes W)=5c_3(TS^6,J)=5e(TS^6)=10v.
\]
Here $H^2(S^6)=H^4(S^6)=0$, and the top Chern class of a complex bundle is its Euler class \cite{MS}. Thus these factors cannot come from $G_{4,5}$, even though $6<2\min(4,5)$. A normal framing and a representative for the label are independent pieces of the cycle data.
\end{example}

\section{Forgetting a framing}\label{sec:forget}
For any CW space $Y$ there are natural homomorphisms
\begin{equation}\label{eq:forget}
 \Omega_d^{\fr}(Y)\xrightarrow{q} A_d(Y)\xrightarrow{h}H_d(Y;\ZZ),
 \qquad h[M,f]=f_*[M].
\end{equation}
The second is well defined because the boundary fundamental class maps to zero in the homology of $W$. For $Y=B$, homology carries the Pontryagin product induced by tensor product of labels. With this convention both maps preserve products.

\begin{theorem}\label{thm:forget}
The map $q$ is onto. Its kernel is generated by differences of compatible-orientation framings on the same labelled manifold. Both maps in \eqref{eq:forget} become isomorphisms over $\QQ$. In particular $\ker q$ is torsion, and $q$ induces an isomorphism of the integral groups modulo their torsion subgroups.
\end{theorem}
\begin{proof}
Every generator of $A_d(Y)$ admits a stable normal framing inducing the given orientation. In the quotient of $\Omega_d^{\fr}(Y)$ by the indicated differences, choosing such a framing gives a well-defined class. For any boundary relation in $A_d(Y)$, choose a stable framing on the oriented bordism. Its restrictions may be used for all boundary components, so the relation also holds in the quotient. This defines an inverse there to the map induced by $q$.

Under Theorem~\ref{thm:thom}, $hq$ is the stable Hurewicz map. Indeed, the ordinary homology Thom isomorphism sends the collapse fundamental class to the fundamental class of $M$, which the final map sends to $f_*[M]$. The rational stable Hurewicz theorem identifies $\pi_d(\Sigma^\infty Y_+)\otimes\QQ$ with $H_d(Y;\QQ)$ \cite[Prop.~10.22]{Cohen}. Thus $(hq)\otimes\QQ$ is an isomorphism. As $q\otimes\QQ$ is surjective, it is an isomorphism and so is $h\otimes\QQ$. Finally, if $q(x)$ is torsion, then some nonzero multiple of $x$ lies in $\ker q$, hence is torsion. It follows that $x$ itself is torsion. This proves the last assertion.
\end{proof}

In particular, the characteristic numbers $\langle f^*z,[M]\rangle$ for $z\in H^d(B;\QQ)$ detect the rational labelled bordism classes in both groups. We next use these numbers with the tensor product multiplication to determine the rational ring.

\section{The rational ring and its sphere representatives}\label{sec:ring}
For $n\ge2$ give $S^{2n}$ its boundary orientation and a stable framing extending over the disk. Choose the generator $b_n\in\widetilde K^0(S^{2n})$ with $\langle\chernchar_n(b_n),[S^{2n}]\rangle=1$, where $\chernchar_j$ has cohomological degree $2j$. Such a normalization follows by expressing a Bott generator as the $n$-fold external product of the reduced Hopf line class on $S^2$, with the sign fixed by the orientation \cite{Atiyah}. Let
\[
 \xi_n=1+b_n,\qquad t_n=[S^{2n},\mathfrak f_{\partial},\xi_n]
           \in\Omega_{2n}^{\fr}(B).
\]
The first Chern class vanishes since $n\ge2$. We also write $t_n$ for $q(t_n)$. On a finite matrix stage, the corresponding factors have classes $k\xi_n$ and $l\xi_n^{-1}$. More generally $\chernchar(f^*u_{k,l})=\chernchar(f^*V_{k,l})/k$ by Lemma~\ref{lem:normalize}.

\begin{theorem}\label{thm:polynomial}
The homomorphisms taking $t_n$ to the sphere classes above give isomorphisms of graded $\QQ$-algebras
\[
 \QQ[t_2,t_3,t_4,\ldots]\ \cong\ \Omega_*^{\fr}(B)\otimes\QQ
       \ \xrightarrow[\cong]{q\otimes\QQ}\ A_*(B)\otimes\QQ,
 \qquad |t_n|=2n.
\]
\end{theorem}
\begin{proof}
As a space $B$ is $\BSU$. Its rational cohomology is polynomial on $\chernchar_2,\chernchar_3,\ldots$, since the Chern character components are related to the usual Chern class generators $c_j$, $j\ge2$, by a triangular change of variables over $\QQ$ \cite{MS,Atiyah}. Thus the monomials indexed by partitions of $n$ with all parts at least two form a basis of $H^{2n}(B;\QQ)$. There is no odd cohomology. By Theorem~\ref{thm:forget} it suffices to exhibit a homology basis for the product induced by tensoring labels.

Let $\lambda=(\lambda_1,\ldots,\lambda_r)$ be such a partition. On the product of the corresponding spheres let $x_i$ be the pullback of the positive top class from factor $i$. Multiplicativity of the Chern character gives
\begin{equation}\label{eq:character}
 \chernchar(\xi_{\lambda_1}\boxtimes\cdots\boxtimes\xi_{\lambda_r})
     =\prod_{i=1}^r(1+x_i),\qquad x_i^2=0.
\end{equation}
In particular $\chernchar_j$ is the sum of the square-free monomials $\prod_{i\in I}x_i$ with $\sum_{i\in I}\lambda_i=j$. For a partition $\mu$ the characteristic number of $\prod_h\chernchar_{\mu_h}$ is the coefficient of $x_1\cdots x_r$. It counts ordered allocations of the labelled factors to blocks having the specified weights $\mu_h$. It is zero unless $\lambda$ refines $\mu$. For $\lambda=\mu$ each block is a singleton and the coefficient is
\[
  \prod_{j\ge2} m_j(\lambda)!,
\]
where $m_j(\lambda)$ is the multiplicity of $j$. The characteristic matrix is therefore triangular with nonzero diagonal after ordering by refinement. Its columns give a rational homology basis. They are precisely the images of the products $\prod_i t_{\lambda_i}$ under $hq$, so those products are a basis in each positive even degree. Degree zero contributes the unit. This proves the ring assertion.
\end{proof}

\begin{example}[Dimension twelve]
For $2n=12$ the characteristic numbers are given by the following table. A column denotes a product of the spheres with their labels, and a row specifies the characteristic class paired with the fundamental class.
\[
\begin{array}{lrrrr}
\toprule
 & (S^4)^3 & (S^6)^2 & S^4\times S^8 & S^{12}\\
\midrule
\chernchar_2^3 & 6 & 0 & 0 & 0\\
\chernchar_3^2 & 0 & 2 & 0 & 0\\
\chernchar_2\chernchar_4 & 3 & 0 & 1 & 0\\
\chernchar_6 & 1 & 1 & 1 & 1\\
\bottomrule
\end{array}
\]
For instance, the entry $3$ counts the choices of the $S^4$ factor contributing to $\chernchar_2$ when the other two contribute jointly to $\chernchar_4$. The determinant is $12$.
\end{example}

\begin{proposition}[Powers of a sphere label]\label{prop:powers}
With the bounding framing, $[S^{2n},\mathfrak f_{\partial},\xi_n^{\,m}]=m t_n$ in $\Omega_{2n}^{\fr}(B)$ for every $m\in\ZZ$. The same identity holds in $A_{2n}(B)$.
\end{proposition}
\begin{proof}
Since $b_n^2=0$, one has $\xi_n^{\,m}=1+m b_n$, the pullback of $\xi_n$ by a based self-map of $S^{2n}$ of degree $m$. Put $d=2n$. The stable basepoint splitting gives $\Omega_d^{\fr}(S^d)=\pi_d^S\oplus\ZZ$, where $\pi_d^S$ denotes the $d$th stable stem. The class of the identity-labelled sphere with its bounding framing is $(0,1)$. Its projection to the point bounds, and its projection to reduced homology is the oriented fundamental class. A based map of degree $m$ acts by multiplication by $m$ on the reduced summand. Pushing forward along the label $S^d\to B$ proves the identity, and applying $q$ gives the second assertion.
\end{proof}

\begin{proposition}[An integral summand in each generator degree]\label{prop:integral}
For $n\ge2$, evaluation of $\chernchar_n$ defines an integer-valued homomorphism on $A_{2n}(B)$ and on $\Omega_{2n}^{\fr}(B)$. In both groups it takes $t_n$ to $1$. Thus $\ZZ t_n$ is an integral direct summand in either group.
\end{proposition}
\begin{proof}
Let $(M,\xi)$ be a cycle of dimension $2n$. Choose an embedding with an oriented trivialized normal bundle of sufficiently large even rank $2r$, and let $c:S^{2n+2r}\to S^{2r}\wedge M_+$ be its collapse. Let $b_r$ be a reduced Bott generator with Chern character the positive cohomology generator $v_r$ of $S^{2r}$. Applying the Chern character to the reduced $K$-class $c^*(b_r\boxtimes\xi)$ gives
\[
 \big\langle\chernchar_{n+r}\big(c^*(b_r\boxtimes\xi)\big),[S^{2n+2r}]\big\rangle
    =\langle\chernchar_n(\xi),[M]\rangle.
\]
Indeed, the degree $2n+2r$ component before pullback is $v_r\boxtimes\chernchar_n(\xi)$, and the ordinary Thom isomorphism evaluates it as indicated. The left side is an integer by Bott normalization on the sphere. The right side depends only on the oriented labelled cycle and vanishes on boundaries. It gives the required homomorphism, and its value on $t_n$ is one by construction. In degree four the same conclusion follows at once from $\chernchar_2=-c_2$ when $c_1=0$.
\end{proof}

The sphere products therefore generate an integral polynomial subring in each theory, since any integral polynomial relation would survive as a relation over $\QQ$. Theorem~\ref{thm:polynomial} identifies the rational rings, and Proposition~\ref{prop:integral} detects each individual sphere integrally. Determining the full integral lattice of products, or the torsion in $A_*(B)$, is a further question.

\paragraph{Use of AI tools.} OpenAI's ChatGPT assisted with reviewing the earlier arguments, developing proofs for this revision, and preparing the manuscript.

\end{document}